\documentclass[12pt]{article}
\usepackage{amssymb}
\DeclareSymbolFont{largesymbol}{OMX}{yhex}{m}{n}
\DeclareMathAccent{\Widehat}{\mathord}{largesymbol}{"62}
\usepackage{amsmath}
\usepackage{tikz}
\usepackage{amscd}
\usepackage{latexsym}
\usepackage{cite}
\DeclareSymbolFont{largesymbol}{OMX}{yhex}{m}{n}
\DeclareMathAccent{\Widehat}{\mathord}{largesymbol}{"62}

\newtheorem{Theorem}{Theorem}[section]

\newtheorem{Lemma}[Theorem]{Lemma}
\newtheorem{Proposition}[Theorem]{Proposition}
\newtheorem{Corollary}[Theorem]{Corollary}

\newtheorem{Question}[Theorem]{Question}

\newcommand{\qed}{\hfill $\Box$}
\newenvironment{Proof}{ \emph{Proof.}}{\qed}
\def\bc{\begin{center}}
\def\ec{\end{center}}

\begin{document}
\title{{\bf The Burnside ai-semiring variety defined by $x^n\approx x$}}
\author{{ \bf Miaomiao Ren$^{1}$}\footnote{E-mail: miaomiaoren@nwu.edu.cn}
   {\bf \quad Xianzhong Zhao $^{1}$\footnote{E-mail: xianzhongzhao@263.net}}
   {\bf \quad Mikhail V. Volkov$^{2}$} \footnote{E-mail: m.v.volkov@urfu.ru}\\
   {\small $^1$ School of Mathematics},
   {\small Northwest University} \\
   {\small Xi'an, Shaanxi, 710127, P.R. China}\\
     {\small $^2$ Institute of Natural Sciences and Mathematics},
   {\small Ural Federal University} \\
   {\small Ekaterinburg, Russia 620000}}
\date{}
\maketitle \vskip -4pt
\baselineskip 16pt
\begin{center}
\emph{To the memory of Professor Libor Pol\'{a}k}
\end{center}
\begin{center}
\begin{minipage}{140mm}
\noindent\textbf{ABSTRACT}
Let ${\bf Sr}(n, 1)$ denote the ai-semiring variety
defined by the identity $x^n\approx x$, where $n>1$.
We characterize all subdirectly irreducible members of a semisimple
subvariety of ${\bf Sr}(n, 1)$. Based on this
result, we prove that ${\bf Sr}(n, 1)$ is hereditarily finitely based (resp., hereditarily finitely generated)
if and only if $n<4$ and that the lattice of subvarieties of ${\bf Sr}(n, 1)$ is countable if and only if $n<4$.
Also, we show that the class of all locally finite members of ${\bf Sr}(n, 1)$
forms a variety and so affirmatively answer the restricted Burnside problem for ${\bf Sr}(n, 1)$.
In addition, we provide a simplified proof of the main result obtained by Gajdo\v{s} and Ku\v{r}il
(Semigroup Forum 80: 92--104, 2010).

\vskip 6pt \noindent
\textbf{Keywords:} Burnside ai-semiring; Burnside group; Variety; Lattice;
Hereditarily finitely based; Hereditarily finitely generated

\vskip 6pt \noindent
\textbf{2010 Mathematics Subject Classifications:}  08B15; 08B26; 16Y60; 20M07
\end{minipage}
\end{center}
\section{Introduction and preliminaries}
By a \emph{semiring} we mean an algebra $(S, +, \cdot)$ of type $(2, 2)$
such that
\begin{itemize}
   \item  the additive reduct $(S, +)$ is a commutative semigroup;
  \item  the multiplicative reduct $(S, \cdot)$ is a semigroup;
  \item $(S, +, \cdot)$ satisfies the identities $x(y+z)\approx xy+xz$, $(y+z)x\approx yx+zx$.
\end{itemize}
One can find many examples of semirings in almost all branches
of mathematics.
Semirings can be regarded as a common generalization
of both rings and distributive lattices. They have
been widely applicated in theoretical computer
science and information science (see \cite{gl, go}).
A semiring is called an \emph{additively idempotent semiring}
(ai-semiring for short, see \cite{rzs1, rzw}) if its additive reduct is a semilattice, i.e.,
a commutative idempotent semigroup.
Such an algebra is also called
a \emph{semilattice-ordered semigroup} (see \cite{gk, kp, rz2}).

Let $X^+$ denote the free semigroup on
a countably infinite set $X$ of variables
and $P_f (X^+)$ the set of all finite non-empty
subsets of $X^+$. If we define
\[
(\forall A,B \in P_f (X^+)) ~ A+B=A\cup B, A\circ B=\{ab \mid a\in A, b\in B\},
\]
then $(P_f (X^+), +, \circ)$ is free in the variety of all ai-semirings
with respect to the mapping $\varphi\colon X \to P_f (X^+), x\mapsto\{x\}$
(see \cite[Theorem 2.5]{kp}).
An ai-semiring identity over $X$ is an
expression of the form $u\approx v$, where $u, v\in P_f(X^+)$.
For convenience, we write
$u_1+\dots+u_k\approx v_1+\dots+v_\ell$ for the ai-semiring identity
$\{u_i \mid 1\leq i\leq k\}\approx \{v_j \mid 1\leq j\leq \ell\}$.

In the last decades, the theory of ai-semiring varieties has been well
developed. In particular, the ai-semiring variety ${\bf Sr}(n, m)$ defined by
the identity $x^n\approx x^m$ with $m<n$ has been intensively studied
(see \cite{gk, gpz, kp, mr, pas1, pz2, rz2, rz3, rzs1, rzs2, rzw, rzs3, zhao1}).
Its members are called \emph{Burnside ai-semirings} (see \cite{rzs1, rzw}).
Recall that a group (resp., semigroup) is called a \emph{Burnside group}
(resp., \emph{Burnside semigroup})
if it satisfies the identity
$x^n\approx x^m$ with $m<n$.
The variety of all Burnside groups
(resp., Burnside semigroups)
satisfying the identity
$x^n\approx x^m$ with $m<n$ will be denoted by
${\bf G}(n, m)$ (resp., ${\bf Sg}(n, m)$).
Therefore, a Burnside ai-semiring is an ai-semiring whose
multiplicative reduct is a Burnside semigroup.

In 1979 McKenzie and  Romanowska \cite{mr}
showed that there are 5 subvarieties of ${\bf Sr}(2, 1)$
satisfying $xy\approx yx$. In 2002
Zhao \cite{zhao1} provided a model of the free object
in ${\bf Sr}(2, 1)$ by using so-called the closed subsemigroups
of the free object in ${\bf Sg}(2, 1)$.
In 2005 Pastijn and Zhao \cite{pz2} showed that
the multiplicative reduct of each member of ${\bf Sr}(2, 1)$
is a regular band, and then Ghosh et al. \cite{gpz} and Pastijn \cite{pas1}
proved that the lattice of subvarieties of ${\bf Sr}(2, 1)$
is a distributive lattice of order 78.
Also, each member of this lattice is finitely based and
finitely generated.

Following this research route
Ku\v{r}il and Pol\'{a}k \cite{kp} initiated to study the variety ${\bf Sr}(n, 1)$.
They provided a construction of the free object in ${\bf Sr}(n, 1)$
by using $n$-closed subsets of the free object in ${\bf Sg}(n, 1)$,
which is a generalization of the notion of the closed subsemigroups.
In 2010 Gajdo\v{s} and Ku\v{r}il \cite{gk}
showed that ${\bf Sr}(n, 1)$ is locally finite if and only if
${\bf G}(n, 1)$ is locally finite. In 2017
Ren et al. \cite{rzs1} showed that the multiplicative reduct of each member of
${\bf Sr}(n, 1)$ is a regular orthocryptogroup. As an application,
models of the free objects in some subvarieties of ${\bf Sr}(n, 1)$ were given
(see \cite{rzs1, rzs2}). Recently, we turn our attention to
the lattice ${\cal L}({\bf Sr}(n, 1))$ of subvarieties of ${\bf Sr}(n, 1)$.
Ren and Zhao \cite{rz2} proved that there are 9 subvarieties of ${\bf Sr}(3, 1)$
satisfying $xy\approx yx$. Based upon this result,
Ren et al. \cite{rzw} showed that the lattice ${\cal L}({\bf Sr}(3, 1))$
is a distributive lattice of order 179.
Also, each member of this lattice is finitely based
and finitely generated.

A variety is said to be \emph{hereditarily finitely based}
(resp., \emph{hereditarily finitely generated})
if all its subvarieties are finitely based (resp., finitely generated). It was proved that
both ${\bf Sr}(2, 1)$ and ${\bf Sr}(3, 1)$ are hereditarily finitely based
and hereditarily finitely generated and that both ${\cal L}({\bf Sr}(2, 1))$ and ${\cal L}({\bf Sr}(3, 1))$
are finite distributive lattices (see \cite{gpz, pas1, rzw}).
Ren et al. \cite{rzs3} showed that the same result holds for
the subvariety
of ${\bf Sr}(n, 1)$ determined by $xy\approx yx$ if $n\geq 4$ and $n-1$ is
square-free.
In this paper we shall continue to study the subvarieties
of ${\bf Sr}(n, 1)$ for $n\geq 4$.
The following questions naturally arise.

\begin{Question}\label{q1}
Is ${\bf Sr}(n, 1)$ hereditarily finitely based for $n\geq 4$$?$
\end{Question}

\begin{Question}\label{q4}
Is ${\bf Sr}(n, 1)$ hereditarily finitely generated for $n\geq 4$$?$
\end{Question}
\begin{Question}\label{q2}
Is the lattice ${\cal L}({\bf Sr}(n, 1))$ finite for $n\geq 4$$?$
\end{Question}

\begin{Question}\label{q3}
Is the lattice ${\cal L}({\bf Sr}(n, 1))$ distributive for $n\geq 4$$?$
\end{Question}

It is well-known that ${\bf Sr}(n, 1)$ is hereditarily finitely based
if and only if ${\cal L}({\bf Sr}(n, 1))$ satisfies the descending chain
condition, i.e., there is no infinite descending chain in ${\cal L}({\bf Sr}(n, 1))$.
Thus an affirmative answer to Question \ref{q2} would imply
an affirmative answer to Question \ref{q1}, since a finite lattice satisfies the descending chain
condition.
On the other hand,
we know that ${\bf Sr}(n, 1)$
is locally finite if and only if ${\bf G}(n, 1)$ is locally finite and
that ${\bf G}(n, 1)$ is not locally
finite for some $n\geq 4$ (see \cite{adi}).
Thus ${\bf Sr}(n, 1)$ is not locally
finite for some $n\geq 4$. This implies that
${\bf Sr}(n, 1)$ is not hereditarily finitely generated for some $n\geq 4$,
since a finitely generated variety is locally finite.
Therefore, to answer the Question \ref{q4}, we need only consider the case where
${\bf Sr}(n, 1)$ is locally finite.
Suppose that ${\bf Sr}(n, 1)$ is locally finite for some $n\geq 4$.
Then it is hereditarily finitely generated if
and only if ${\cal L}({\bf Sr}(n, 1))$ satisfies the ascending chain
condition, i.e., there is no infinite ascending chain in ${\cal L}({\bf Sr}(n, 1))$
(see \cite{mv}). Hence an affirmative answer to Question \ref{q2} would imply
an affirmative answer to Question \ref{q4}.

Let ${\bf V}$ be a variety. The so-called \emph{Burnside problem} for ${\bf V}$
means whether it is locally finite.
Novikov and Adian \cite{adi} gave a negative answer to the Burnside problem for the group variety ${\bf G}(n, 1)$
for all even $n>665$.
Also, the Burnside problem for the semigroup variety ${\bf Sg}(n, 1)$
has been reduced to a group-theory problem, due to Green and Rees \cite{gr}:
${\bf Sg}(n, 1)$ is locally finite if and only if ${\bf G}(n, 1)$ is locally finite.
Based upon this result, Gajdo\v{s} and Ku\v{r}il \cite{gk} showed that
the Burnside problem for ${\bf Sr}(n, 1)$ is equivalent to
that for ${\bf G}(n, 1)$.
Our interest is to simplify its proof.

Magnus \cite{mag} asked whether all locally finite members of ${\bf G}(n, 1)$ forms a
group variety  for any positive integer $n$.
Based on the classification of finite simple groups, Hall and Higman \cite{hh} reduced
this problem to the case when $n=p^k +1$ ($p$ is a prime)
and then Zel'manov \cite{zel1, zel2} completely solved it.
Later on, Sapir asked whether all locally finite members of a variety ${\bf V}$ forms
a variety (see \cite[Problem 3.10.3]{sa3}) and called it the
\emph{restricted Burnside problem} for ${\bf V}$. He \cite{sa2} provided an algorithm
to decide whether the restricted Burnside problem for a finitely based semigroup
variety has positive answer.
In this paper, we shall solve  the restricted Burnside problem for ${\bf Sr}(n, 1)$.

In addition to this section, this paper is organized as follows.
In Sect. 2 we shall characterize the subdirectly irreducible members of
the subvariety ${\bf M}_n$ of ${\bf Sr}(n, 1)$ and show that it is semisimple.
In Sect. 3 we shall prove that the class of all locally finite members of ${\bf Sr}(n, 1)$
forms a variety and so affirmatively answer the restricted Burnside problem for ${\bf Sr}(n, 1)$.
Also, we shall prove that ${\bf Sr}(n, 1)$ is hereditarily finitely based (resp., hereditarily finitely generated)
if and only if $n<4$ and that the lattice of subvarieties of ${\bf Sr}(n, 1)$ is countable if and only if $n<4$.
Moreover, we shall provide a simplified proof of the main result obtained by Gajdo\v{s} and Ku\v{r}il \cite{gk}.

For other notations and terminology used in this paper, the reader
is referred to Petrich and Reilly \cite{pet} for
a background on semigroup theory, to Ghosh et al.
\cite{gpz},
Pastijn \cite{pas1}, Pastijn and Zhao \cite{pz2} and Ren et al. \cite{rzw}
for knowledge on ai-semiring varieties, and to Burris and Sankappanavar \cite{bur} for
information concerning universal algebra. We shall assume that the
reader is familiar with the basic results in these areas.

\section{A semisimple subvariety of ${\bf Sr}(n, 1)$}
Let ${\bf M}_n$ denote the subvariety of ${\bf Sr}(n, 1)$ determined by the identity
\begin{align}\label{new3}
x^{n-1}+y^{n-1}\approx x^{n-1}y^{n-1}.
\end{align}
In this section we shall characterize the subdirectly irreducible members of ${\bf M}_n$
and prove that it is semisimple, i.e., all its subdirectly irreducible members are congruence
simple. The following lemma will be useful for us later.
\begin{Lemma}\label{lemm1.13}
${\bf M}_n$ satisfies the following identities:
\begin{align}
x^{n-1}y^{n-1}  & \approx y^{n-1}x^{n-1};    \label{new5}\\
xy^{n-1}        & \approx y^{n-1}x;    \label{228}\\
x+y             & \approx xy^{n-1}+x^{n-1}y;\label{new1}\\
x+x^{n-1}       & \approx x+x^2+\cdots+x^{n-1};\label{21701}\\
(xy)^{n-1}       & \approx x^{n-1}y^{n-1}.\label{new2}
\end{align}
\end{Lemma}
\begin{Proof}
Suppose that $S$ is a member of ${\bf M}_n$.

For any $a, b \in S$, we have
\[a^{n-1}b^{n-1} \stackrel{(\ref{new3})} = a^{n-1}+b^{n-1} = b^{n-1}+a^{n-1} \stackrel{(\ref{new3})} = b^{n-1}a^{n-1}.\]
This shows that ${\bf M}_n$ satisfies the identity (\ref{new5}).

As usual, we denote Green's
${\cal D}$-relation and ${\cal H}$-relation on
$(S, \cdot)$ by ${\cal D}$ and ${\cal H}$, respectively. Then
for any $a, b \in S$, by \cite[Lemma 4.1]{rzs1} we have
\[
(a, b)\in \mathcal{D}\Leftrightarrow a^{n-1}b^{n-1}a^{n-1}=a^{n-1},b^{n-1}a^{n-1}b^{n-1}=b^{n-1}
\]
and
\[
(a, b)\in \mathcal{H}\Leftrightarrow a^{n-1}=b^{n-1}.
\]
This implies that  ${\cal D}$ is equal to ${\cal H}$,
since ${\bf M}_n$ satisfies the identity (\ref{new5}).
Also, it is easy to see that $(S, \cdot)$ is completely regular. Therefore, by \cite[Theorem IV.2.4]{pet} we have that $(S, \cdot)$ is a Clifford semigroup and so ${\bf M}_n$ satisfies the identity (\ref{228}).

Suppose that $a, b \in S$. Then we have
\begin{align*}
a+b
&= (a+b)^n \\
&= (a+b)^n + ab^{n-1}+a^{n-1}b \\
&= a+b+ab^{n-1}+a^{n-1}b\\
&= a^n+b^n+ab^{n-1}+a^{n-1}b\\
&= a(a^{n-1}+b^{n-1})+(a^{n-1}+b^{n-1})b\\
&= aa^{n-1}b^{n-1}+a^{n-1}b^{n-1}b&&(\text{by}~(\ref{new3})) \\
&= ab^{n-1}+a^{n-1}b
\end{align*}
and
\begin{align*}
&a+a^{n-1}\\
&= (a+a^{n-1})^n\\
&= (a+a^{n-1})^n+a(a^{n-1})^{n-1}+a^2(a^{n-1})^{n-2}+\cdots+a^{n-1}a^{n-1}\\
&= (a+a^{n-1})^n+a+a^2+\cdots+a^{n-1}\\
&= (a+a^{n-1})+a+a^2+\cdots+a^{n-1}\\
&= a+a^2+\cdots+a^{n-1}.
\end{align*}
This shows that ${\bf M}_n$ satisfies the identities (\ref{new1}) and (\ref{21701}).

From \cite{rzs1} we know that ${\bf Sr}(n, 1)$ satisfies the identity (\ref{new2}).
It follows that so does for ${\bf M}_n$, since it is a subvariety of ${\bf Sr}(n, 1)$.
\end{Proof}

By a 0-group $G^0$ we mean the semigroup obtained from a group $G$ by adding an extra
zero element 0. Let $S$ be a member of ${\bf M}_n$ and $E(S)$ denote the set of all idempotents of $(S, \cdot)$.
Then by \cite[Lemma 2.1]{rzs1} $E(S)$ is a subsemiring of $S$.

\begin{Lemma}\label{lemm1.15} Let $S$ be a member of ${\bf M}_n$.
If the multiplicative reduct of $S$ is a $0$-group, then $S$ is congruence simple.
\end{Lemma}
\begin{Proof}
Suppose that the multiplicative reduct of $S$ is a 0-group $G^0$, where
$1_G$ is the identity of the group $G$. Then $E(S)=\{1_G, 0\}$.
Let $a$ be an element of $G\backslash\{1_G\}$. Then it is easy to verify that $a+a^2+\cdots+a^{n-1}\in E(S)$.
If $a+a^2+\cdots+a^{n-1}=1_G$, then
\[a=a1_G=a(a+a^2+\cdots+a^{n-1})=a+a^2+\cdots+a^{n-1}=1_G,\]
a contradiction. This implies that $a+a^2+\cdots+a^{n-1}=0$, since $E(S)=\{1_G, 0\}$.
By (\ref{21701}) it follows that $a+a^{n-1}=0$. Thus
\begin{align}\label{102}
a\in G \backslash\{1_G\} \Rightarrow a+a^{n-1}=0.
\end{align}
Let $c$ and $d$ be different elements of $S$. It is easy to check
that $cd^{n-2}\neq 1_G$.
Then we have
\begin{align*}
c+d
&= cd^{n-1}+c^{n-1}d&&(\text{by}~(\ref{new1}))\\
&= (cd^{n-2}+c^{n-1}d^{n-1})d\\
&= (cd^{n-2}+(cd^{n-2})^{n-1})d&&(\text{by}~(\ref{new2}))\\
&= 0d&&(\text{by}~(\ref{102}))\\
&=0.
\end{align*}
That is to say,
\begin{align}\label{216}
c \neq d \Rightarrow c+d=0.
\end{align}
By \cite[Theorem 4.7]{mf} it follows that $S$ is congruence simple.
\end{Proof}
\begin{Lemma}\label{lemm1.16}
Let $S$ be a member of ${\bf M}_n$. Then a semiring congruence $\theta$ on $S$ is
diagonal if and only if $\theta|_{E(S)}$ is the identity congruence on $E(S)$.
\end{Lemma}
\begin{Proof}
Suppose that $\theta$ is the identity congruence on $S$. Then it is easy to see that
$\theta|_{E(S)}$ is the identity congruence on $E(S)$.

Conversely, assume that $\theta$ is a semiring congruence on $S$ such that
$\theta|_{E(S)}$ is diagonal on $E(S)$. Let $a$ be an arbitrary element of $S$.
If $(a, a^{n-1})\in \theta$, then $(a^{n-1}+a, a^{n-1})\in \theta$.
Also,  by (\ref{21701}) it follows that both $a^{n-1}+a$ and $a^{n-1}$
are in $E(S)$.
This implies that $a^{n-1}+a=a^{n-1}$. Again by (\ref{21701}) we have
\begin{align*}
a=aa^{n-1} & =a(a^{n-1}+a)=a(a+a^2+\cdots+a^{n-1}) \\
           & =a+a^2+\cdots+a^{n-1} \\
           & =a^{n-1}+a \\
           & =a^{n-1}.
\end{align*}
Now we have proved that $(a, a^{n-1})\in \theta \Rightarrow a=a^{n-1}$. This shows that $\theta|_{G_\alpha}$ is diagonal for all $\alpha \in Y$,
where the multiplicative reduct of $S$ is the Clifford semigroup $[G_\alpha, Y]$.
On the other hand, if $(a, b)\in \theta$ for some $a, b\in S$,
then $(a^{n-1}, \, b^{n-1})\in \theta$ and so $a^{n-1}=b^{n-1}$, since $a^{n-1}, \, b^{n-1}\in E(S)$ and $\theta|_{E(S)}$ is diagonal on $E(S)$.
This implies that $a, b \in G_\alpha$ for some  $\alpha \in Y$ and so $a=b$.
We therefore have proved that $\theta$ is the identity congruence.
\end{Proof}

The following lemma tells us that each  semiring congruence on
$E(S)$ can be extended to a congruence on $S$.
\begin{Lemma}\label{lemm1.17}
Let $S$ be a member of ${\bf M}_n$. If $\rho$ is a semiring congruence on
$E(S)$, then there exists a congruence $\tau$ on
$S$ such that the restriction $\tau|_{E(S)}=\rho$.
\end{Lemma}
\begin{Proof}
Define a binary relation $\tau$ on $S$ by
\[
(a, b)\in \tau  \Leftrightarrow (\exists~ e\in E(S))~ ea=eb,~(e,  a^{n-1})\in \rho, (a^{n-1}, b^{n-1})\in \rho.
\]
By (\ref{new5}) and (\ref{228})
it is easy to check that $\tau$ is a semigroup congruence on
$(S, \cdot)$ and that the restriction $\tau|_{E(S)}=\rho$. In the remainder
we shall show that $\tau$ is a semigroup congruence on $(S, +)$.
We first prove that
\begin{align}\label{002}
(a, f)\in \tau \textrm{ for some} ~f \in E(S) \Rightarrow (a, a+a^{n-1})\in\tau.
\end{align}
In fact, if $(a, f)\in \tau$ for some $f \in E(S)$, then there exists
$e\in E(S)$ such that $ea=ef$, $(e, a^{n-1})\in \rho$ and  $(a^{n-1}, f)\in \rho$.
This implies that $(a, ef)\in \tau$
and so $(a(a+a^{n-1}), ef(a+a^{n-1}))\in \tau$.
By (\ref{new5}) and (\ref{21701}) we have that $ef(a+a^{n-1})=ef$
and that $a(a+a^{n-1})=a+a^{n-1}$.
Thus $(a+a^{n-1}, ef)\in\tau $ and so $(a, a+a^{n-1})\in\tau$.

Next, we shall show that
\begin{align}\label{003}
(a, b)\in\tau\Rightarrow (a+a^{n-1}, b+b^{n-1})\in \tau.
\end{align}
Assume that $(a, b)\in\tau$. Then $(a(b+b^{n-1}), b(b+b^{n-1}))\in\tau$. Also,
by (\ref{21701}) it follows that $b(b+b^{n-1})=b+b^{n-1}$ and so
$(a(b+b^{n-1}), b+b^{n-1})\in\tau$.
Since $b+b^{n-1} \in E(S)$,
by (\ref{002}) we have
\[
(a(b+b^{n-1}), a(b+b^{n-1})+(a(b+b^{n-1}))^{n-1})\in\tau.
\]
Also, by (\ref{new2}) we have
\begin{align*}
a(b+b^{n-1})+(a(b+b^{n-1}))^{n-1}
& =a(b+b^{n-1})+a^{n-1}(b+b^{n-1})^{n-1}\\
& =a(b+b^{n-1})+a^{n-1}(b+b^{n-1})\\
& =(a+a^{n-1})(b+b^{n-1}).
\end{align*}
It now follows that
\[
(a(b+b^{n-1}), (a+a^{n-1})(b+b^{n-1})))\in\tau.
\]
This implies that
\[
(b+b^{n-1}, (a+a^{n-1})(b+b^{n-1})))\in\tau.
\]
Similarly,
\[
(a+a^{n-1}, (b+b^{n-1})(a+a^{n-1})))\in\tau.
\]
Since $a+a^{n-1}, b+b^{n-1} \in E(S)$,
by (\ref{new5}) it follows that
\[(a+a^{n-1})(b+b^{n-1})=(b+b^{n-1})(a+a^{n-1})\]
and so $(a+a^{n-1}, b+b^{n-1})\in \tau$.

Finally, suppose that $(a, b)\in\tau$. Then for any $c \in S$,
$(ac^{n-2}, bc^{n-2})\in \tau$. By (\ref{003}) we have
\[
(ac^{n-2}+(ac^{n-2})^{n-1}, bc^{n-2}+(bc^{n-2})^{n-1})\in \tau.
\]
This implies that
\[
((ac^{n-2}+(ac^{n-2})^{n-1})c, (bc^{n-2}+(bc^{n-2})^{n-1})c)\in \tau.
\]
Notice that, by (\ref{new1}) and (\ref{new2}),
\begin{align*}
(ac^{n-2}+(ac^{n-2})^{n-1})c
& =(ac^{n-2}+a^{n-1}c^{n-1})c\\
& =ac^{n-1}+a^{n-1}c=a+c
\end{align*}
and
\begin{align*}
(bc^{n-2}+(bc^{n-2})^{n-1})c
& =(bc^{n-2}+b^{n-1}c^{n-1})c\\
& =bc^{n-1}+b^{n-1}c=b+c.
\end{align*}
Hence
$(a+c, b+c)\in\tau$.
\end{Proof}

Suppose that $S$ is a member of ${\bf M}_n$. By (\ref{new5}) we can define
the partial order relations $\leq_+$ and $\leq_\cdot$ on $S$ as follows:
\begin{align*}
a\leq_+ b    &\Leftrightarrow a+b=b,\\
a\leq_\cdot b&\Leftrightarrow (\exists~ e \in E(S)) ~a=eb.
\end{align*}
If $a\leq_\cdot b$, then $a=eb$ for some $e\in E(S)$. Further,
\[a=eb=e^{n-1}b^{n-1}b \stackrel{(\ref{new2})}=(eb)^{n-1}b=a^{n-1}b,\]
i.e., $a=a^{n-1}b$. We therefore have
\[a+b=a^{n-1}b+b^{n-1}b=(a^{n-1}+b^{n-1})b\stackrel{(\ref{new3})}=a^{n-1}b^{n-1}b=a^{n-1}b=a.\]
This shows that $b\leq_+ a$. Conversely, if $a\leq_+ b$ on $S$, then $b=a+b$ and so
$b^{n-1}=(a+b)^{n-1}$. Further, we have
\begin{align*}
b^{n-1}
&= (a+b)^{n-1}= (a+b)^{n-1}+ab^{n-2}+a^{n-1}\\
&=b^{n-1}+ab^{n-2}+a^{n-1}
\end{align*}
and
\begin{align*}
b^{n-1}
&= (a+b)^{n-1}= (a+b)^{n-1}+a^{n-2}b+b^{n-2}a\\
&=b^{n-1}+a^{n-2}b+b^{n-2}a
\end{align*}
That is to say,
\begin{align}
b^{n-1}    &= b^{n-1}+ab^{n-2}+a^{n-1}\label{new10}
\end{align}
and
\begin{align}
b^{n-1}    &= b^{n-1}+a^{n-2}b+b^{n-2}a  \label{new11}
\end{align}
It follows that
\begin{align*}
b=a+b
&=(a+b)^{n-1}(a+b)=b^{n-1}(a+b)\\
&=(b^{n-1}+ab^{n-2}+a^{n-1})(a+b)&&(\text{by}~(\ref{new10}))\\
&=(a^{n-1}b^{n-1}+ab^{n-2})(a+b)&&(\text{by}~(\ref{new3}))\\
&=ab^{n-1}+a^{n-1}b+ab^{n-2}a   &&(\text{by}~(\ref{228}))\\
&=a(b^{n-1}+a^{n-2}b+b^{n-2}a)\\
&=ab^{n-1}.&&(\text{by}~(\ref{new11}))
\end{align*}
Thus $b\leq_\cdot a$ and so $\leq_+=_\cdot\geq$. To summarize, we
have proved
\begin{Lemma}\label{newlemma}
Let $S$ be a member of ${\bf M}_n$. Then
$\leq_+=_\cdot\geq$ holds for $S$.
\end{Lemma}

\begin{Lemma}\label{lemm1.18}
If $S$ be a subdirectly irreducible member of
${\bf M}_n$, then the multiplicative reduct of $S$ is a
$0$-group.
\end{Lemma}
\begin{Proof}
Suppose that $\theta$ is the least non-identity semiring congruence on
$S$. Then by Lemmas \ref{lemm1.16} and \ref{lemm1.17},
$\theta|_{E(S)}$ is the least non-identity semiring congruence on
$E(S)$. This implies that $E(S)$ is a subdirectly
irreducible member of $\mathbf{M}_2$ and so it has exactly two elements.
Suppose now that $E(S)=\{0,\, 1 \}$ with $0 \leq_\cdot 1$ and that
$a$ is an arbitrary element of $S$. By Lemma \ref{newlemma} it is easily verified
that the greatest lower bound $ a \wedge 1$ of $\{a, 1\}$ under $\leq_\cdot$
exists. Then there
exist $e,\,f\in E(S)$ such that $a\wedge 1= e a$ and $a\wedge
1=1f=f$. This implies that $ea=f$. Thus
\[0=0 f=0(ea)=(0 e)a=0a.\]
By (\ref{228}) it follows that $0$ is the
zero element of $(S, \cdot)$.
Since $(S, \cdot)$ is completely regular,
we have that the
${\cal H}$-classes of $(S, \cdot)$ containing $0$ and $1$ are
$\{0\}$ and $S\backslash\{0\}$, respectively. Thus $(S, \cdot)$ is
a 0-group as required.
\end{Proof}

By Lemmas \ref{lemm1.15} and \ref{lemm1.18} we now have

\begin{Proposition}\label{pro1}
Let $S$ be a member of ${\bf M}_n$. Then the following statements are equivalent:
\begin{itemize}
  \item [$(a)$]     $S$ is subdirectly irreducible;
  \item [$(b)$]    $S$ is congruence simple;
  \item [$(c)$]   the multiplicative reduct of $S$ is a $0$-group.
\end{itemize}
\end{Proposition}

As a consequence, we have
\begin{Corollary}
${\bf M}_n$ is semisimple.
\end{Corollary}

Let $(G, \cdot)$ be a member of ${\bf G}(n, 1)$. Then $(G^0, \cdot)$
becomes a member of ${\bf Sg}(n, 1)$. Define an addition on $G^0$ by
the rule
\begin{align}\label{201}
a+b=
\begin{cases}
0& if ~a\neq b, \\
a& if ~a=b.
\end{cases}
\end{align}
It is easy to verify that $(G^0, +, \cdot)$ is a member of ${\bf M}_n$.
The algebra $G^0$ is called the \emph{flat extension} of $G$ or a \emph{flat semiring}
(see \cite{jac}).
By Proposition \ref{pro1} and (\ref{216}) we immediately deduce
\begin{Corollary}\label{npro1}
If $G$ is a group in ${\bf G}(n, 1)$, then the flat extension of $G$ is
a subdirectly irreducible member of ${\bf M}_n$.

Conversely, every subdirectly irreducible member of
${\bf M}_n$ is isomorphic to the flat extension of some member of ${\bf G}(n, 1)$.
\end{Corollary}

\section{Main Results}
Let ${\cal L}_q({\bf G}(n, 1))$ denote the lattice of subquasivarieties of ${\bf G}(n, 1)$
and ${\bf HSP}({\cal A})$ the variety generated by a class ${\cal A}$ of ai-semirings.
Since a subclass of ${\bf G}(n, 1)$ is a universal Horn class
if and only if it is a quasivariety, by \cite[Corollary 5.4]{jac} and Corollary \ref{npro1}
we have
\begin{Proposition}\label{pro01}
Let $\psi$ be a mapping defined by the rule
\[
\psi: {\cal L}_q({\bf G}(n, 1)) \to [{\bf M}_2, {\bf M}_n],
{\cal V}\mapsto {\bf HSP}(\{G^0\mid G\in {\cal V}\}).
\]
Then $\psi$ is a lattice isomorphism.
\end{Proposition}

We immediately obtain
\begin{Corollary}\label{co01}
The lattice ${\cal L}_q({\bf G}(n, 1))$ can be embedded into ${\cal L}({\bf Sr}(n, 1))$.
\end{Corollary}

By \cite[Theorem 7.3]{jac} we have
\begin{Lemma}\label{lem01}
Let $G$ be a finite member of ${\bf G}(n, 1)$. Then
the flat semiring $G^0$ is finitely based if and only if
all Sylow subgroups of $G$ are abelian.
\end{Lemma}

For each odd prime $p$ one can consider the $p^3$-element group
\[
G_p=\langle a,b,c \mid a^p=b^p=c^p=e, ab=ba,\ bc=cb,\ cac^{p-1}=ab \rangle,
\]
where $e$ is the identity of $G_p$.
Each element of $G_p$ admits a unique representation
of the form $a^ib^jc^k$ where $0\le i,j,k<p$ and the
multiplication in $G_p$ may be expressed by the formula
\[
a^ib^jc^k\cdot a^mb^nc^r=
a^{i+m(\bmod p)}b^{j+km+n(\bmod p)}c^{k+r(\bmod p)}.
\]
It is known that $G_p$ is a non-abelian group in ${\bf G}(p+1, 1)$.
By Lemma \ref{lem01} we have that
the flat semiring $G_p^0$ is a nonfinitely based member of ${\bf Sr}(p+1, 1)$.
This implies that ${\bf Sr}(p+1, 1)$ is not hereditarily finitely based.
On the other hand, consider the quaternion group
\[
Q_8=\langle i,j,k \mid i^{2}=j^{2}=k^{2}=ijk=e\rangle,
\]
where $e$ is the identity of $Q_8$.
It is easy to see that $Q_8$ is a 8-element non-abelian group
in ${\bf G}(5, 1)$.
By Lemma \ref{lem01} we have that the flat semirings $Q^0_8$
is a nonfinitely based member of ${\bf Sr}(5, 1)$.
This implies that ${\bf Sr}(5, 1)$ is not hereditarily finitely based.

Suppose now that $n\geq4$. Then either $n-1$ has an odd prime divisor $p$
or it is divisible by 4. It follows that either ${\bf Sr}(n,1)$ contains
${\bf Sr}(p+1,1)$ for some odd prime number $p$ or it contains
${\bf Sr}(5,1)$. This implies that
${\bf Sr}(n, 1)$ is not hereditarily finitely based.
Notice that both ${\bf Sr}(2, 1)$ and ${\bf Sr}(3, 1)$
are hereditarily finitely based (see \cite{gpz, pas1, rzw}).
We have proved the following theorem:
\begin{Theorem}\label{theo2}
${\bf Sr}(n, 1)$ is hereditarily finitely based if and only if $n<4$.
\end{Theorem}

The above theorem negatively answers Question \ref{q1}. This enables us to
provide a negative answer to Question \ref{q2}. Indeed, we have

\begin{Theorem}\label{theo1}
The lattice ${\cal L}({\bf Sr}(n, 1))$ is countable if and only if $n<4$.
\end{Theorem}
\begin{Proof}
If $n<4$, then ${\cal L}({\bf Sr}(n, 1))$ is a finite lattice (see \cite{gpz, pas1, rzw}).
Suppose now that $n\geq4$. Then either ${\bf Sr}(n, 1)$ contains ${\bf Sr}(p+1, 1)$ for some
odd prime number $p$ or ${\bf Sr}(n, 1)$ contains ${\bf Sr}(5, 1)$.
This implies that either ${\cal L}({\bf Sr}(n, 1))$ has a sublattice ${\cal L}({\bf Sr}(p+1, 1))$
for some odd prime number $p$
or it has the sublattice ${\cal L}({\bf Sr}(5, 1))$.
Therefore, to show that  ${\cal L}({\bf Sr}(n, 1))$ is uncountable,
it is enough to prove that
${\cal L}({\bf Sr}(5, 1))$ and ${\cal L}({\bf Sr}(p+1, 1))$
are uncountable for all odd prime number $p$. In fact,
by the main theorem of \cite{fed} we have that for each odd prime $p$, there exists a finite
member of ${\bf G}(p+1, 1)$ which
that generates a quasivariety with a continuum number of
subquasivarieties. This shows that the lattice ${\cal L}_q({\bf G}(p+1, 1))$
is uncountable. By Corollary \ref{co01} it follows that
the lattice ${\cal L}({\bf Sr}(p+1, 1))$ also is uncountable.
On the other hand, if a locally finite
group variety has countably many subquasivarieties, then it
must consist of groups whose nilpotent subgroups are all abelian
(see \cite[Theorem 1]{sa}). Since the group
variety ${\bf G}(5,1)$ is locally finite (see \cite{san}) and
contains non-abelian 8-element groups (the quaternion group,
say), it follows that the lattice ${\cal L}_q({\bf G}(5, 1))$ is
uncountable. Again by Corollary \ref{co01} we have that
${\cal L}({\bf Sr}(5, 1))$ also is uncountable.
This completes the proof.
\end{Proof}

\begin{Theorem} \label{theo3}
${\bf Sr}(n, 1)$ is hereditarily finitely generated if and only if $n<4$.
\end{Theorem}
\begin{Proof}
If $n<4$, then ${\bf Sr}(n, 1)$ is hereditarily finitely generated (see \cite{pas1, rzw}).
Suppose now that $n\geq4$. Consider the following two cases:

\textbf{Case 1.} ${\bf Sr}(n, 1)$ is not locally finite. Then
it is not hereditarily finitely generated, since a finitely generated variety must be locally finite.

\textbf{Case 2.} ${\bf Sr}(n, 1)$ is locally finite.
Then it is hereditarily finitely generated if
and only if ${\cal L}({\bf Sr}(n, 1))$ satisfies the ascending chain
condition (see \cite{mv}). Thus, to prove that ${\bf Sr}(n, 1)$ is not hereditarily finitely generated,
it is enough to show that ${\cal L}({\bf Sr}(n, 1))$ does not satisfy the ascending chain
condition. We know that either ${\cal L}({\bf Sr}(n, 1))$ has a
sublattice ${\cal L}({\bf Sr}(p+1, 1))$ for some odd prime number $p$
or it has the sublattice ${\cal L}({\bf Sr}(5, 1))$.
Therefore, it suffices to prove that ${\cal L}({\bf Sr}(5, 1))$ and ${\cal L}({\bf Sr}(p+1, 1))$
does not satisfy the ascending chain condition for all odd prime number $p$.
From \cite{quick} we have that the lattice of subvarieties of ${\bf G}(5, 1)$ contains an infinite
ascending chain and so does ${\cal L}_q({\bf G}(5, 1))$. By Corollary \ref{co01} it follows that
${\cal L}({\bf Sr}(5, 1))$ does not satisfy the ascending chain condition.
On the other hand,
for every odd prime number $p$, ${\cal L}({\bf G}(p+1, 1))$ contains an infinite
ascending chain (see \cite{kn}) and so does ${\cal L}_q({\bf G}(p+1, 1))$.
Again by Corollary \ref{co01} we have that ${\cal L}({\bf Sr}(p+1, 1))$ does not
satisfy the ascending chain condition.

This completes the proof.
\end{Proof}

Theorem \ref{theo3} negatively answers Question \ref{q4}. Next, we shall
provide a simplified proof of \cite[Theorem 2.4]{gk}. For this, the following auxiliary result
is necessary.

\begin{Lemma}\label{nle1}
An ai-semiring is locally finite if and only if its multiplicative reduct is locally finite.
\end{Lemma}
\begin{Proof}
Suppose that $A$ is a finite subset of an ai-semiring $S$. We shall use $\langle A \rangle$ and
$\langle A \rangle_s$ to denote the subsemiring of $S$ generated by $A$
and the subsemigroup of $(S, \cdot)$ generated by $A$,
respectively.  Then it is easy to verify that
\[
\langle A \rangle=\{b_1+b_2+\cdots b_m \mid b_i \in \langle A \rangle_s,
i=1,2,\ldots,m, m\geq 1\}
\]
and so $\langle A \rangle$ is finite if and only if $\langle A \rangle_s$
is finite. This shows that $S$ is locally finite if and only if $(S, \cdot)$ is locally finite.
\end{Proof}

\begin{Proposition}
The following statements are equivalent:
\begin{itemize}
  \item [$(a)$]    ${\bf Sr}(n, 1)$ is locally finite;
  \item [$(b)$]   ${\bf Sg}(n, 1)$ is locally finite;
  \item [$(c)$]  ${\bf G}(n, 1)$ is locally finite.
\end{itemize}
\end{Proposition}
\begin{Proof}
$(a)\Rightarrow(c)$. Suppose that ${\bf G}(n, 1)$ is not locally finite.
Then there exists an infinite finitely generated group $G$ in ${\bf G}(n, 1)$
and so the semigroup $(G^0, \cdot)$ is not locally finite.
By Lemma \ref{nle1} we have that
the flat semiring $G^0$ is not locally finite and so is ${\bf Sr}(n, 1)$.
This contradicts the hypothesis that ${\bf Sr}(n, 1)$ is locally finite.
Thus ${\bf G}(n, 1)$ is locally finite.

$(c)\Rightarrow(b)$. This follows from the main result of \cite{gr}.

$(b)\Rightarrow(a)$. Let $S$ be an arbitrary member of ${\bf Sr}(n, 1)$. Then
$(S, \cdot)$ is a semigroup in ${\bf Sg}(n, 1)$. By $(b)$ we have that
$(S, \cdot)$ is locally finite. It follows from Lemma \ref{nle1}
that $S$ is locally finite and so is ${\bf Sr}(n, 1)$.
\end{Proof}

In the remainder of this section we shall affirmatively answer the restricted Burnside problem for ${\bf Sr}(n, 1)$.
The following result, which is due to Sapir \cite[Theorem 1]{sa2}, is necessary.
\begin{Lemma}\label{lem124}
Let {\bf V} be a periodic variety of semigroups that is finitely based. Then
there is a recursive function $f(k)$ bounding the orders of $k$-generator finite
semigroups in {\bf V} if and only if all
nil-semigroups in ${\bf V}$ are locally finite.
\end{Lemma}

By \cite[Problem 3.10.4]{sa3} we have that the class of
all locally finite members of a finitely based variety {\bf V}
forms a variety
if and only if for every $k\geq 1$, there is an upper bound for the size of all $k$-generated
finite algebras in {\bf V}. It is obvious that ${\bf Sr}(n, 1)$ is finitely based.
Thus, to prove that the class of all locally finite members of ${\bf Sr}(n, 1)$
forms a variety, it is enough to show that for every $k\geq 1$, there is an
upper bound for the size of all $k$-generated
finite semirings in ${\bf Sr}(n, 1)$.
Suppose now that $\langle A \rangle$ is a
$k$-generated finite semirings in ${\bf Sr}(n, 1)$. Then
$\langle A \rangle_s$ is a $k$-generated finite semigroup in ${\bf Sg}(n, 1)$
and by the proof of Lemma \ref{nle1} we have
\begin{align}\label{n222}
|\langle A \rangle|\leq 2^{|\langle A \rangle_s|}.
\end{align}
Since ${\bf Sg}(n, 1)$ is a periodic variety and
all nil-semigroups in this variety are trivial semigroups,
it follows from Lemma \ref{lem124} that there is an upper bound $f(k)$ for the size of all $k$-generated
finite semigroups in ${\bf Sg}(n, 1)$ and so $|\langle A \rangle_s|\leq f(k)$.
By (\ref{n222}) we deduce that $|\langle A \rangle|\leq 2^{f(k)}$.
Thus we have

\begin{Theorem}
The class of all locally finite members of ${\bf Sr}(n, 1)$
forms a variety.
\end{Theorem}
{\bf{Conclusion}}
We have negatively answered Questions \ref{q1}, \ref{q4} and \ref{q2}. However,
Question \ref{q3} is still open. We conjecture that the lattice ${\cal L}({\bf Sr}(n, 1))$
is not distributive for all $n\geq4$.
By Corollary \ref{co01} and the proof of
Theorem \ref{theo1} it is enough to show that
${\cal L}_q({\bf G}(5, 1))$ and ${\cal L}_q({\bf G}(p+1, 1))$
are not distributive for all odd prime number $p$.
By \cite[Theorem 5]{bg} we know that the lattice of subquasivarieties
of a quasivariety is distributive if and only if it is modular.
Thus, to prove that ${\cal L}({\bf Sr}(n, 1))$ is not distributive for all $n\geq4$,
it remains to show that ${\cal L}_q({\bf G}(5, 1))$ and ${\cal L}_q({\bf G}(p+1, 1))$
are not modular for all odd prime number $p$.

Theorems \ref{theo2}, \ref{theo1} and \ref{theo3} make us realize the difficulty of
studying subvarieties of ${\bf Sr}(n, 1)$ for $n\geq 4$.
In the next work,
we shall continue to study the related problems of these subvarieties.
The ultimate aim is to classify them with respect to the finitely based property
and the finitely generated property, respectively.

{\bf Acknowledgments.}
Miaomiao Ren is supported by National Natural Science
Foundation of China (11701449) and
Natural Science Foundation of Shaanxi Province (2020JM-009).
Xianzhong Zhao is supported by National Natural Science
Foundation of China (11971383, 11571278). Mikhail Volkov is supported by the Russian Science Foundation (grant No. 22-21-00650)


\begin{thebibliography}{aa}
\bibitem{adi} S.I. Adian, \emph{The Burnside Problem and Identities in Groups}
(Springer, Berlin, 1979).

\bibitem{bg}
A.I. Budkin and V.A. Gorbunov, Quasivarieties of algebraic systems,
\emph{Algebra and Logic}, {\bf 14}(2) (1975), 73--84.

\bibitem{bur} S. Burris and H.P. Sankappanavar, \emph{A Course in Universal Algebra}
(Springer Verlag, New York, 1981).

\bibitem{fed} A.N. Fedorov, Subquasivarieties
of nilpotent minimal non-Abelian group varieties,
\emph{Siberian Math. J.}, {\bf 21}(6) (1980), 840--850.

\bibitem{gk}  P. Gajdo\v{s} and M. Ku\v{r}il,
On free semilattice-ordered semigroups satisfying $x^n\approx x$,
\emph{Semigroup Forum}, {\bf 80} (2010), 92--104.

\bibitem{gpz}  S. Ghosh, F. Pastijn and X.Z. Zhao,
Varieties generated by ordered bands I,
\emph{Order}, {\bf 22} (2005), 109--128.

\bibitem{gl}  K. G{\l}azek, \emph{A Guide to the Literature on Semirings
and their Applications in Mathematics and Information Science}
(Kluwer Academic Publishers, Dordrecht-Boston-London, 2001).

\bibitem{go}  J.S. Golan, \emph{The Theory of Semirings with Applications
in Mathematics and Theoretical Computer Science}
(Longman Scientific and Technical, Harlow, 1992).

\bibitem{gr}  J.A. Green and D. Rees, On semigroups in which $x^r=x$,
\emph{Proc. Camb. Philos. Soc.}, {\bf 48} (1952), 35--40.

\bibitem{hh} P. Hall and G. Higman, On the $p$-length of $p$-soluble groups and reduction
theorems for Burnside's problem, \emph{Proc. London Math. Soc.,} {\bf 6}(3) (1956), 1--42.

\bibitem{jac} M. Jackson, Flat algebras and the translation of universal Horn logic to equational logic,
\emph{J. Symbolic Logic}, {\bf 73}(1) (2008), 90--128.

\bibitem{kn} L.G. Kov\'{a}cs and M.F. Newman,
On non-cross varieties of groups, \emph{J. Austral. Math. Soc.}, {\bf 12}(2) (1971), 129--144.

\bibitem{kp} M. Ku\v{r}il and L. Pol\'{a}k,
On varieties of semilattice-ordered semigroups, \emph{Semigroup Forum}, {\bf 71} (2005), 27--48.

\bibitem{mf}
S.S. Mitchell and P.B. Fenoglio, Congruence-free commutative semirings,
\emph{Semigroup Forum}, {\bf 37} (1988), 79--91.

\bibitem{mv} S.O. Macdonald and M.R. Vaughan-Lee, Varieties that make one Cross,
\emph{J. Austral. Math. Soc. (Series A)}, {\bf 26} (1978), 368--382.

\bibitem{mag}
W. Magnus, A collection between the Baker-Hausdorff formula and a problem of
Burnside,
\emph{Ann. Math.}, {\bf 52} (1950), 11--26; Errata \emph{Ann. Math.} {\bf 57} (1953), 606.

\bibitem{mr} R. McKenzie and A. Romanowska, Varieties of $\cdot$-distributive
bisemilattices, \emph{Contrib. Gen. Algebra}
(Proc. Klagenfurt Conf. 1978) {\bf 1} (1979), 213--218.

\bibitem{pas1} F. Pastijn, Varieties generated by ordered bands II,
\emph{Order}, {\bf 22} (2005), 129--143.

\bibitem{pz2} F. Pastijn and X.Z. Zhao,
Varieties of idempotent semirings with commutative addition, \emph{Algebra Universalis},
{\bf 54} (2005), 301--321.

\bibitem{pet} M. Petrich and N.R. Reilly, \emph{Completely Regular Semigroups} (Wiley, New
York, 1999).

\bibitem{quick} M. Quick,
A classification of some insoluble varieties of groups of exponent four,
\emph{J. Algebra},
{\bf 197}(2) (1997), 342--371.

\bibitem{rz2} M.M. Ren and X.Z. Zhao, The varieties of semilattice-ordered semigroups satisfying
$x^3\approx x$ and $xy\approx yx$, \emph{Period. Math. Hungar.},
{\bf 72} (2016), 158--170.

\bibitem{rz3} M.M. Ren and X.Z. Zhao, Some results on varieties of additively idempotent semirings,
\emph{Pure and Applied Mathematics},
{\bf 34}(4) (2018), 406--410. (In chinese)

\bibitem{rzs1} M.M. Ren, X.Z. Zhao and Y. Shao,
On a variety of Burnside ai-semirings satisfying $x^{n}\approx x$,
\emph{Semigroup Forum}, {\bf 93}(3) (2016), 501--515.

\bibitem{rzs2}
M.M. Ren, X.Z. Zhao and Y. Shao, Varieties of Burnside ai-semirings satisfying $x^{n}\approx x$,
\emph{Acta Sci. Math. (Szeged)} {\bf 87}(1--2) (2021), 39--61.

\bibitem{rzw} M.M. Ren, X.Z. Zhao and A.F. Wang,
On the varieties of ai-semirings satisfying $x^{3}\approx x$,
\emph{Algebra Universalis}, {\bf 77} (2017), 395--408.

\bibitem{rzs3} M.M. Ren, X.Z. Zhao and Y. Shao,
The lattice of ai-semiring varieties satisfying $x^n\approx x$ and $xy \approx yx$
\emph{Semigroup Forum}, {\bf 100} (2020), 542--567.

\bibitem{san}  I.N. Sanov, Solution of the Burnside problem for exponent 4,
\emph{U\v{c}en. Zap. Leningrad Univ.}, {\bf 10} (1940), 166--170.

\bibitem{sa1}  M.V. Sapir, On cross semigroup varieties and related questions,
\emph{Semigroup Forum}, {\bf 42}(1) (1984), 345--364.

\bibitem{sa}  M.V. Sapir, Varieties with a countable number of subquasivarieties,
\emph{Siberian Math. J.}, {\bf 25}(3) (1984), 461--473.

\bibitem{sa2}  M.V. Sapir, The restricted Burnside problem for varieties of semigroups,
\emph{Izv. Akad. Nauk SSSR Ser. Math.}, {\bf 55}(3) (1991), 670--679.

\bibitem{sa3} M.V. Sapir, \emph{Combinatorial Algebra: Syntax and Semantics}, (Springer Monogr. Math., 2014).

\bibitem{zel1} E.I. Zel'manov, The solution of the restricted Burnside problem for groups of odd
exponent, \emph{Izv. Akad. Nauk. SSSR. Ser. Mat.}, {\bf 54}(1) (1990), 42--59.

\bibitem{zel2} E.I. Zel'manov,
The solution of the restricted Burnside problem for $2$-groups,
\emph{Mat. Sb. (N.S.)}, {\bf 182}(4) (1991), 568--592.

\bibitem{zhao1} X.Z. Zhao, Idempotent semirings with a commutative additive reduct,
\emph{Semigroup Forum}, {\bf 64} (2002), 289--296.
\end{thebibliography}
\end{document}